# GEOMETRIES FOR GROUPS

ROBERT ALONZO LYMAN

Abstract. We prove a "continuous converse" to the Milnor-Schwarz lemma. Given any topological group $G$ and any left-invariant, large scale geodesic pseudometric $d$ giving an identity neighborhood finite diameter (for example, the word metric with respect to an open generating set), we construct a geodesic metric space $X$ equipped with a continuous, isometric, cobounded action of $G$ such that $(G, d)$ and $X$ are quasi-isometric via any orbit map.

## 1. Introduction

Suppose a subset $S \subset G$ generates the group $G$. Every element $g \in G$ may be written as a *word* $g = s_1 \cdots s_n$ in the elements of $S$ and their inverses, and the left-invariant *word metric* $d_S : G \times G \to \mathbb{R}$ assigns to a pair of elements $g$ and $h$ in $G$ the minimal word length of $g^{-1}h$.

Discussing this situation in his monograph "Asymptotic Invariants of Infinite Groups" [Gro93], Gromov writes,

> The distance function [...] makes $G$ subject to a geometric scrutiny as any other metric space.
>
> This space may appear boring and uneventful to a geometer's eye since it is discrete and the traditional local (e.g. topological and infinitesimal) machinery does not run in $G$. To regain the geometric perspective, one has to change his/her position and move the observation point far away from $G$. Then the metric in $G$ seen from the distance $d$ becomes the original distance divided by $d$ and for $d \to \infty$ the points in $G$ coalesce into a connected continuous solid unity which occupies the visual horizon without any gaps of holes and fills our geometer's heart with joy.
>
> — Mikhail Gromov [Gro93]

### 1.1. The Milnor-Schwarz Lemma..

Here, Gromov likely has in mind the construction of the asymptotic cone of a metric space, but, discovering geometric group theory as an undergraduate, I found this metaphor a compelling way to think about the relation of quasi-isometry, and noted in particular the appearance of word metrics, which crop up when discussing groups and quasi-isometry because of the Milnor-Schwarz Lemma. We reproduce here a classical version of the statement.

**Lemma 1.1.** (Schwarz [Šva55], Milnor [Mil68]) *Suppose a group $G$ acts properly discontinuously, compactly and by isometries on a proper geodesic metric space $(X, d)$. Then $G$ is finitely generated and if $S$ is any finite generating set, the group $G$ equipped with the word metric $d_S$ and the space $(X, d)$ are quasi-isometric via any orbit map.*

Although classical properness more generally and proper discontinuity play essentially no role in this paper, let us note quickly that this definitions (despite the name) is really a statement about *continuous* actions of topological groups; in this case that the map $G \times X \to X \times X$ sending $(g, x)$ to $(x, g.x)$ is a proper map when $G$ has the discrete topology.

The classic Milnor-Schwarz lemma also has a converse: supposing that $G$ is finitely generated, say by a finite set $S$, then there exists a geodesic metric space $(X, d)$ on which $G$ acts by isometries such that $(G, d_S)$ and $(X, d)$ are quasi-isometric via any orbit map. One such space $X$ is the *Cayley graph* with respect to the (symmetrized) generating set $S \cup S^{-1}$.

The purpose of this note is to discuss the situation of left-invariant, continuous pseudometrics on a general topological group $G$. We prove the following theorem.

**Theorem A.** *Suppose that $\partial$ is a left-invariant, large-scale geodesic pseudometric on a group $G$, for example a word metric $d_S$ for some generating set $S$. There is a geodesic metric space $(X, d)$ with a cobounded, isometric action of $G$ such that $(G, \partial)$ and $(X, d)$ are quasi-isometric via any orbit map.*

*If there is an identity neighborhood $U \ni 1$ in $G$ with finite $\partial$-diameter (for example if the generating set $S$ is an identity neighborhood), then the action of $G$ on $X$ may be taken to be continuous.*

A pseudometric $\partial$ on a space $X$ is *large-scale geodesic* if there exists $R \geq 1$ such that for any two points $x$ and $y$ in $X$, we can find a sequence $x = x_0, ..., x_n = y$ such that successive distances satisfy $d(x_i, x_{i+1}) \leq R$ and we have

$$\sum_{i=1}^{n} d(x_{i-1}, x_i) \leq Rd(x, y).$$

One point of the theorem for the author is similar to the reaction of Gromov's hypothetical geometer: thinking merely of word metrics feels (to me) somewhat "boring and uneventful" without the promise of a possible geodesic metric space to realize them.

Also, since the author trained in the theory of finitely generated groups, the topological situation occasionally feels vexed, if not dangerous. For example, as abelian groups, $\mathbb{R}$, $\mathbb{C}$, $\mathbb{R}^n$ for each $n > 1$, and $\mathbb{Q}_p$ are isomorphic, so only retain their geometric differences when considered as topological groups acting *continuously.* It is fortunate that the condition allowing our construction to run is not at all technical.

`Theorem A` may be compared to realization results in the world of Polish groups. In [GK03], Gao and Kechris show that every Polish group may be realized as the isometry group of some Polish metric space. Melleray [Mel08] showed that compact metrizable groups are isometry groups of compact metric spaces, and Malicki and Solecki [MS09] showed that locally compact, second countable groups are isometry groups of proper Polish metric spaces.

Like many of these results, our proof of `Theorem A` uses the space of *Katětov functions* $E(X)$ on a metric space $X$. These are functions $f : X \to \mathbb{R}_{\geq 0}$ satisfying the inequalities

$$|f(x) - f(y)| \leq d(x, y) \leq f(x) + f(y)$$

for all pairs of points $x$ and $y$ in $X$. The space $E(X)$ is (affinely isomorphic to) a subspace of $\mathrm{BC}(X)$, the Banach space of bounded continuous functions on $X$ equipped with the $\ell^\infty$ supremum norm. Since $\mathrm{BC}(X)$, like all normed spaces, is geodesic, it is therefore quite possible (although not immediately obvious) that the metric spaces constructed in the papers above are already geodesic. On the other hand, although our construction to prove `Theorem A` is *more* flexible, since we can realize any (topologically relevant) word metric, we make no attempt to restrict the group of isometries of the resulting space.

For a reader with a background in finitely generated groups, let us remark that the difficulty in the theorem is that in the general topological setting, constructions involving Cayley graphs are unavailable: the vertex set of a metric graph is discrete, so the stabilizer of a vertex in a nontrivial continuous action is a proper open subgroup of $G$. When $G$ is connected, for instance, no such subgroup exists.

Most of the remainder of the note is given to our proof of `Theorem A`. A few consequences for the Rosendal theory of large-scale geometry of groups are drawn in a final section.

1.2. **Acknowledgments.**

The author is grateful to Carlos Adrian Perez Estrada, Christian Rosendal, Eduardo Martínez-Pedroza, Michael Kopreski, Brian Udall and Jeremy Brazas for helpful conversations on ideas and questions surrounding the content of this note, as well as for a careful reading and helpful feedback on a previous version.

The statement and general methods of proof of `Lemma 2.8` were suggested by an LLM, ChatGPT Pro. The lemma should not be surprising to experts, but comes from an area the author is still building familiarity with. None of the text of this paper is LLM-generated.

## 2. Katětov functions of finite support

Suppose $X$ is a metric space. A function $f : X \to \mathbb{R}_{\geq 0}$ is *Katětov* if it satisfies the following inequalities for each $x$ and $y$ in $X$.

$$|f(x) - f(y)| \leq d(x,y) \leq f(x) + f(y).$$

The space $E(X)$ of Katětov functions on $X$ was studied by Katětov [Kat88] to construct Urysohn universal homogeneous metric spaces. It plays a role in [Fim+22], [Mel08], and [MS09].

Here are some basic properties. The following lemmas are standard, use only the definition and the triangle inequality and are left to the reader.

**Lemma 2.1.** *Given $x \in X$, the function $\iota_x(y) = d(x,y)$ is Katětov, and letting $x$ vary defines a map $\iota : X \to E(X)$ defined as $\iota(x) = \iota_x$ which is an isometric embedding with respect to the $\ell^\infty$ metric*

$$\|f - g\|_\infty = \sup\{|f(x) - g(x)| : x \in X\}.$$

■

If $G$ acts on $X$ on the left by isometries, then $G$ acts on $E(X)$ by the rule

$$g.f(x) = f(g^{-1}.x).$$

In particular

$$g.\iota_x(y) = \iota_x(g^{-1}.y) = d(x, g^{-1}.y) = d(g.x, y) = \iota_{g.x}.$$

This action is by isometries of the $\ell^\infty$ metric. When $G$ acts continuously on $X$, the action on $E(X)$ is also continuous, since for each $f$, letting $g \to 1$ we have

$$\|g.f - f\|_\infty = \sup_{x \in X} |f(g^{-1}.x) - f(x)| \leq \sup_{x \in X} d(g^{-1}.x, x),$$

which tends to 0 since $G$ acts continuously.

**Lemma 2.2.** *Suppose $f \in E(X)$ is Katětov. For each $x \in X$, we have*

$$\|f - \iota_x\|_\infty = f(x).$$

■

Let $S \subset X$ and consider functions $\omega : S \to \mathbb{R}_{\geq 0}$. For each choice of $S$ and $\omega$, we define a function $f = f(S, \omega) : X \to \mathbb{R}_{\geq 0}$ as follows.

$$f(x) = \inf\{\omega(s) + d(x, s) : s \in S\}.$$

We call $S$ the *support* of $f$ and the numbers $\omega(s)$ the *weights* of $f$.

**Lemma 2.3.** *Given $f = f(S, \omega)$ as above, the function $f$ is Katětov provided that for each pair of points $s$ and $t \in S$, we have $d(s, t) \leq \omega(s) + \omega(t)$.*

We remark that as defined above, the same function $f$ may be represented with different choices of supports and weights. We require merely the existence of supports and weights in this paper.

The following is a consequence of `Lemma 2.2`.

**Corollary 2.4.** *Given $f = f(S, \omega)$ as above,*

$$\inf\{\|\iota_x - f\|_\infty : x \in X\} = \inf\{\omega(s) : s \in S\}.$$

*If the righthand infimum is a* minimum *realized by $s \in S$, then the left-hand infimum is realized by $\|\iota_s - f\|_\infty = f(s) = \omega(s)$.* ∎

*Proof.* Suppose $f(x) < f(y)$ and $f(x) \leq \omega(s) + d(x, s)$. Then $f(y) \leq \omega(s) + d(y, s)$, so

$$|f(x) - f(y)| = f(y) - f(x) \leq |d(s, y) - d(s, x)| \leq d(x, y).$$

On the other hand, for each $s$ and $t$ in $X$, by our assumption on $\omega(s) + \omega(t)$, we have

$$d(x, y) \leq d(x, s) + d(s, t) + d(t, y) \leq (\omega(s) + d(s, x)) + (\omega(t) + d(t, y)).$$

Taking infimums completes the proof. □

**Lemma 2.5.** *Fix $S \subset X$ and a closed interval $I \subset \mathbb{R}$ and suppose that for each $s \in S$, we have a 1-Lipschitz function $t \mapsto \omega_t(s)$ so that the function $f_t(x) = \inf\{\omega_t(s) + d(x, s) : s \in S\}$ is Katětov. The path $I \to E(X)$ defined as $t \mapsto f_t$ is 1-Lipschitz, i.e.*

$$\|f_t - f_u\|_\infty \leq |t - u|$$

*for $t$ and $u$ in $I$.*

*Proof.* Suppose first that the support $S = \{s\}$ in the statement, that is, that for each $t$, the function $f_t$ is $f_t(x) = \omega(t) + d(x, s)$. Then for all $x \in X$ we have

$$|f_t(x) - f_u(x)| = |\omega(t) - \omega(u)| \leq |t - u|,$$

since $\omega_t(s) = \omega(t)$ is 1-Lipschitz by assumption. Since the inequality holds for general $S$ after taking infimums, the result follows. □

**Definition 2.6.** For $\varepsilon > 0$, let $X^+ = X^+(\varepsilon)$ be the space of Katětov functions $f(S, \omega)$ of finite support (meaning $S$ is a finite set) with minimum weight

$$\min\{\omega(s) : s \in S\} < \varepsilon.$$

We equip $X^+(\varepsilon)$ with the length metric that assigns to a pair of points $x$ and $y$ in $X^+(\varepsilon)$ the infimal length (measured with the $\ell^\infty$ distance) of a path connecting them.

The following is an immediate consequence of the definitions.

**Lemma 2.7.** *If for each pair of points $x$ and $y$ in $X$, there is a sequence*

$$x = x_0, ..., x_n = y$$

*such that for $i$ satisfying $1 \leq i \leq n$ we have $d(x_{i-1}, x_i) < \varepsilon$, then $X^+ = X^+(\varepsilon)$ is a length space satisfying $d_X(x, y) \leq d_{X^+}(\iota_x, \iota_y)$, and $\iota(X)$ is $\varepsilon$-dense in $X^+$.* ■

**Lemma 2.8.** *If $X$ is a length space, then for all $\varepsilon > 0$ we have that $X^+(\varepsilon)$ is a length space. Furthermore, the embedding $\iota : X \to X^+$ is isometric, and for each pair of points $x$ and $y$ in $X$ there is a geodesic between $\iota_x$ and $\iota_y$ in $X^+$.*

*Proof.* Fix $0 < \delta < \frac{\varepsilon}{2}$. Since $X$ is a length space, for each $x$ and $y$ in $X$, there exists a sequence $x = x_0, ..., x_n = y$ such that for $i$ satisfying $1 \leq i \leq n$ we have $d(x_{i-1}, x_i) < \delta$, and such that

$$\sum_{i=1}^{n} d(x_{i-1}, x_i) < d(x, y) + \delta.$$

Let $S = \{x_0, ..., x_n\}$, let $D = d(x, y)$, and let $S_k$ be the partial sum $\sum_{i=1}^{k} d(x_{i-1}, x_i)$. For $t \in [0, D]$ define $\omega_t : S \to \mathbb{R}_{\geq 0}$ as

$$\omega_t(s) = \begin{cases} t & s = x_0 = x \\ \delta + |t - S_i| & s = x_i, \text{ for } \ 1 \leq i \leq n-1 \\ D - t & s = x_n = y. \end{cases}$$

The reader can easily verify that the pairwise sums $\omega_t(x_i) + \omega_t(x_j)$ satisfy the hypotheses of `Lemma 2.3`, so for each $t \in [0, D]$, the function $f_t = f(S, \omega_t)$ is Katětov. By `Lemma 2.5`, we have

$$\|f_t - f_u\|_\infty \leq |t - u|$$

for each $t$ and $u$ in $[0, D]$.

The reverse inequality holds by considering $|f_t(x) - f_u(x)|$, so the curve $t \mapsto f_t$ is a geodesic of length $D = d(x, y)$.

Evaluating $\|f_0 - \iota_x\|_\infty = f_0(x)$ and $\|f_D - \iota_y\|_\infty = f_D(y)$ shows that $f_t$ starts at $x$ and ends at $y$. Finally, for each $t \in [0, D]$, there exists $i$ satisfying $0 \leq i \leq n$ such that $\omega_t(x_i) < 2\delta < \varepsilon$, so this geodesic from $x$ to $y$ of length $d(x, y)$ is contained in $X^+$, completing the proof. □

Before turning to the proof of `Theorem A`, we need the following variant of [Ros22, Lemma 2.70].

**Lemma 2.9.** *Suppose that $\partial$ is a left-invariant, large-scale geodesic pseudometric on a group $G$ and that there exists an identity neighborhood with finite $\partial$-diameter. There exists a* continuous *left-invariant pseudometric $\hat{\partial}$ on $G$ quasi-isometric to $\partial$.*

*Proof.* Suppose that $\partial$ is $R$-large-scale geodesic. The standard proof of the Milnor-Schwarz lemma (see e.g. [BH99]) proves that $\partial$ is quasi-isometric to the word metric $d_S$ with respect to the symmetric generating set $S = \{g \in G : \partial(1, g) < R\}$. Let $U$ be a symmetric identity neighborhood with finite $d_S$-diameter. The word metric with respect to the *open,* symmetric generating set $US \cup SU$ is bi-Lipschitz equivalent to $d_S$ by the triangle inequality and the assumption that $U$ has finite $d_S$-diameter.

There exists a *continuous,* left-invariant pseudometric $d$ assigning $US \cup SU$ finite diameter, and the proof of [Ros22, Lemma 2.70] shows that

$$\hat{\partial}(g, h) = \inf\left\{ \sum_{i=1}^{n} d(x_{i-1}, x_i) : g = x_0, ..., x_n = h \right\}$$

is a continuous, left-invariant pseudometric quasi-isometric to the word metric with respect to $US \cup SU$ and hence to $\partial$. □

Finally, we use `Lemma 2.8` to complete the proof of `Theorem A`.

*Proof of* **`Theorem A`**. First, use `Lemma 2.9` to replace the given left-invariant, large-scale geodesic pseudometric $\partial$ with a continuous left-invariant pseudometric $\hat{\partial}$. Since $\partial$ was assumed to be large-scale geodesic, so is $\hat{\partial}$. Let

$$H = \left\{ h \in G : \hat{\partial}(1, h) = 0 \right\}.$$

This is a closed subgroup of $G$, and the quotient space $Y = G/H$ equipped with the quotient metric induced by $\hat{\partial}$ is a large-scale geodesic metric space equipped with a continuous, transitive action of $G$.

Since $Y$ is large-scale geodesic, let $R > 0$ such that for all $x$ and $y$ in $Y$, there exists a sequence $x = x_0, ..., x_n = y$ of points in $Y$ such that for $1 \leq i \leq n$ we have $d(x_{i-1}, x_i) < R$ and

$$\sum_{i=1}^{n} d(x_{i-1}, x_i) < Rd(x, y).$$

By `Lemma 2.7`, the space $Y^+(R)$ is a length space containing $\iota(Y)$ as an $R$-dense subset, and the existence of the sequence above implies that the metric on $Y^+(R)$ satisfies

$$d_Y(x, y) \leq d_{Y^+}\left(\iota_x, \iota_y\right) \leq Rd_Y(x, y).$$

The continuous, transitive, isometric action of $G$ on $Y$ yields a continuous, cobounded, isometric action on the length space $Y^+(R)$.

Now we restore geodesics to $Y^+(R)$ by repeatedly applying `Lemma 2.8`. Let $X_0 = Y^+(R)$ and for each $n \geq 1$, let $X_n = (X_{n-1})^+\left(\frac{1}{2^n}\right)$ and consider the direct limit

$$X_\infty = \lim_{n \to \infty} X_n.$$

By `Lemma 2.8`, the embedding $\iota : X_n \to X_{n+1}$ is isometric, and each pair of points of $\iota(X_n)$ are joined by a geodesic in $X_{n+1}$.

The isometric embedding of $X_0$ in $X_n$ for each $n$ yields an isometric embedding into $X_\infty$, and each point of $X_\infty$ is at distance at most $\sum_n \frac{1}{2^n}$ from $X_0$, hence $R + \sum_n \frac{1}{2^n}$ from the image of $Y$. By the above argument, the continuous, cobounded, isometric action of $G$ on $X_0$ extends to $X_\infty$. Each pair of points of $X_\infty$ belong to some $X_n$ and are joined by a geodesic in $X_{n+1}$, so the metric space $X_\infty$ is geodesic and still satisfies

$$d_Y(x, y) \leq d_{X_\infty}\left(\iota_x, \iota_y\right) \leq Rd_Y(x, y).$$

□

## 3. Consequences for coarse geometry of topological groups

Let us begin this section by recalling some generalities. A *topological group* is a group $G$ equipped with a topology for which multiplication $m : G \times G \to G$ and inversion $(\cdot)^{-1} : G \to G$ are continuous maps ($G \times G$ has the product topology).

**Lemma 3.1.** (Pontryagin [Pon39]) *The topological closure of the identity in $G$ is a closed normal subgroup $N \lhd G$ and the quotient $G/N$ is a Hausdorff group. Every topological group $G$ fits into a short exact sequence*

$$1 \to C \to G \to T \to 1$$

*where each homomorphism is continuous, $C$ is closed in $G$ and $T$ has the quotient topology.*

*The group $C$ is connected (and contains the subgroup $N$ above), and the group $T$ has a neighborhood basis of the identity given by open subgroups.*

Such groups $T$ are often called *non-Archimedean* or *zero-dimensional* and are totally disconnected. In the case where $G$ is non-Archimedean, `Theorem A` is much easier: let $U$ be an open subgroup with finite $d_S$-diameter, and let $\Gamma$ be the simplicial graph with vertex set $G/U$, and an edge from $gU$ to each distinct vertex of the form $gusU$ for $u \in U$ and $s \in S$.

The action of $G$ on $\Gamma$ is continuous, simplicial and vertex transitive, hence isometric when $\Gamma$ has the geodesic graph metric. It's not difficult to show that

$$d_\Gamma(gU, hU) \leq d_S(g, h) \leq (M + 1)d_\Gamma(gU, hU) + M,$$

where $M = \sup_{u \in U} d_S(1, u)$.

In [Ros22], Rosendal studies the canonical *left-coarse structure* on a topological group $G$. This is the coarse structure generated by the collection of sets which are bounded with respect to some continuous, left-invariant pseudometric on $G$. A subset $B$ of $G$ is said to be *coarsely bounded* if it is bounded in this coarse structure. This means that in every continuous action of $G$ by isometries on a metric space (or equivalently in every continuous, left-invariant pseudometric), the subset $B$ has bounded orbits (equivalently bounded diameter).

Of particular interest in this setting are the groups which are *locally bounded,* which is to say have a coarsely bounded identity neighborhood, and those which are *monogenic,* being algebraically generated by a coarsely bounded set. An immediate consequence of `Theorem A` is the following.

**Corollary 3.2.** (Continuous Milnor–Schwarz) *A topological group $G$ is monogenic and locally bounded if and only if there exists a geodesic metric space $X$ on which $G$ acts continuously, coboundedly, and* metrically properly discontinuously *by isometries. In this situation, $G$ has a well-defined quasi-isometry type, which is realized by the action on $X$ via any orbit map.*

The condition that $G$ acts *metrically properly discontinuously* says that for each $x \in X$ and $R > 0$, the set

$$\{g \in G : g.B_R(x) \cap B_R(x) \neq \varnothing\}$$

is coarsely bounded in $G$.

*Proof.* The "if" direction is proven in [Ros22] and follows the classical proof (say, in [BH99]) more or less verbatim. The "only if" direction follows from `Theorem A` by constructing the geodesic metric space $X$ associated to an *open, coarsely bounded* generating set $S$, which exists by assumption. □

Let us mention a little of the relative setting (used classically when studying finitely generated relatively hyperbolic groups, for example). Suppose that $G$ is a topological group and that $\mathcal{H}$ is a finite collection of closed subgroups. We say that $G$ is *monogenic relative to $\mathcal{H}$* if there is a symmetric, coarsely bounded set $S$ such that $G$ is generated by $T = S \cup \bigcup \mathcal{H}$. If the word metric $d_T$ gives an identity neighborhood finite diameter (for example, if $G$ is locally bounded and $S$ contains an identity neighborhood), then `Theorem A` provides a geodesic metric space $X$ on which $G$ acts coboundedly, continuously and isometrically such that $(X, d)$ is quasi-isometric to $(G, d_T)$ via any orbit map.

Arora and Martínez-Pedroza [AM23] study the case where $G$ is a topological group with a compact open subgroup and $\mathcal{H}$ is comprised of open subgroups.

Consider the simplicial graph $\Gamma$ with vertex set $G \cup \bigcup_{h\in\mathcal{H}} G/H$ and edges connecting a vertex $g$ to all vertices of the form $gs$ for $s \in S$ and to all vertices $gH$ for $H \in \mathcal{H}$. This graph $\Gamma$ is quasi-isometric to $(G, d_T)$. In the discrete setting, the graph $\Gamma$ is called a *coned-off Cayley graph* [Far98, Hru10].

In the totally disconnected locally compact or non-Archimedean Polish settings, a related construction when $\mathcal{H}$ is the empty collection is called a *Cayley-Abels graph* [ÁLM23] or a *rough Cayley graph* [KM08] in the TDLC case and a *Cayley-Abels-Rosendal graph* by Branman, Domat, Hoganson and the present author in the non-Archimedean Polish setting [Bra+25].

In all of the above settings, the graph itself (or a graph quasi-isometric to it) suffices, since $G$ will act on the geodesic metric space $\Gamma$ continuously. In the case where $\mathcal{H}$ is a single closed $G_\delta$ subgroup $H$, a clever use of the inputs into the proof of `Theorem A` can produce a geodesic metric space $X$ where stabilizers are finite intersections of conjugates of $H$. It is interesting to note that this is a far cry from the specificity with which Arora and Martínez-Pedroza are able to understand the action of $G$ on their graph $\Gamma$.

Here are some questions.

**Question 3.3.** Suppose $\Gamma = \Gamma(S, \mathcal{H})$ as described above where $H$ comprises at least two closed subgroups. Can the proof of `Theorem A` produce a geodesic metric space $X$ where stabilizers are either coarsely bounded in $G$ or a finite intersection of conjugates of subgroups in $\mathcal{H}$?

**Question 3.4.** The action of $G$ on the geodesic metric space $X$ constructed in `Theorem A`, although cobounded, is never transitive. Yet many connected locally compact groups (for example the classical semi-simple Lie groups) act continuously, isometrically and transitively on geodesic metric spaces. To what extent can the same be said for connected groups in general?

**Question 3.5.** Likewise, the geodesic metric space $X$ constructed in `Theorem A` is likely never uniquely geodesic. Which topological groups act continuously, coboundedly and isometrically on *uniquely* geodesic metric spaces?

**Question 3.6.** Suppose that some of the spaces $X$ constructed in `Theorem A` for the simplicial graph $\Gamma = \Gamma(S, \mathcal{H})$ described above are $\delta$-hyperbolic in the sense of Gromov. It is well-known that an additional property, such as fineness [Bow12] or bounded coset penetration [Far98] is necessary for a study of relative hyperbolicity. One challenge in adopting the definition verbatim in our setting is that $\Gamma$ is frequently uncountable, even of uncountable valence. Can fineness or bounded coset penetration be studied in $\Gamma$ directly? Are there analogous properties for $X$ which suffice?

Department of Mathematics and Computer Science, Rutgers University–Newark, Newark, NJ 07102
*Email address:* robbie.lyman@rutgers.edu
*URL:* https://robbielyman.com/math